\documentclass[12pt,leqno]{amsart}

\usepackage{amssymb}
\usepackage{latexsym}

\newcommand\dspace{\lineskip=2pt\baselineskip=18pt\lineskiplimit=0pt}
\dspace

\newcounter{deficislo}[section]
\newcounter{propcislo}[section]
\newcounter{lemacislo}[section]
\newcommand{\defi}[1]{\refstepcounter{deficislo}{\noindent \bf
Definition \thesection.\thedeficislo.\  }{\it #1}}
\newcommand{\pro}[1]{\refstepcounter{propcislo}{\noindent \bf
Proposition \thesection.\thepropcislo.\  }{\it #1}}
\newcommand{\lema}[1]{\refstepcounter{lemacislo}{\noindent \bf
Lemma \thesection.\thelemacislo.\  }{\it #1}}

\newcommand{ \qz}{\tilde P(z,\bar z, u)}

\newcommand{ \tg}{\tilde g}

\newcommand{\de}{\delta}
\newcommand{\ab}[1]{\vert z\vert^{#1}}
\newcommand{\cdva}{{\mathbb C^2}}

\newcommand{\te}{\theta}

\newcommand{\fz}{{F(z,\bar z,u)}}

\newcommand{ \al}{\alpha}

\newcommand{\tf}{\tilde f}

\begin{document}
\title{
 Local equivalence of  symmetric hypersurfaces in
 $\mathbb C^2$.}
\author{Martin Kol\'a\v r}
\address {Department of Mathematics and Statistics, Masaryk University,
Janackovo nam. 2a, 662 95 Brno } \email {mkolar@math.muni.cz }

\maketitle

\begin{abstract} The Chern-Moser normal form
and its analog on  finite type hypersurfaces in general do not
respect symmetries. Extending the work of N. K. Stanton, we
consider the local equivalence problem for symmetric Levi
degenerate hypersurfaces of finite type in $\mathbb C^2$. The
results give for all such hypersurfaces a complete normalization
which respects the symmetries. In particular, they apply to tubes
and rigid hypersurfaces, providing an effective classification.
The main tool is a complete normal form constructed
 for a general hypersurface with a tube model.
As an application, we describe all biholomorphic maps between
tubes, answering a question posed by N. Hanges.
 Similar results for hypersurfaces admitting
nontransversal  symmetries are obtained.
\end{abstract}

\section{Introduction}

One of fundamental problems in CR geometry concerns local
biholomorphic classification of real hypersurfaces in complex
space.
 An extrinsic approach to the problem, originating in the work
 of Poincar\'e
(\cite{Po}),  is to  analyze directly the action of local
biholomorphisms on the defining equation of the hypersurface.

 In
the Levi nondegenerate case this approach was completed in the
normal form construction of  Chern and  Moser (\cite{CM}). A
similar construction  for Levi degenerate hypersurfaces of finite
type in dimension two was obtained in (\cite{Ko1}).

 As an
immediate application, normal forms can be used
 for a closely related
geometric problem, to   determine  local symmetries of a
hypersurface. In fact, except for the sphere and its blow-ups, all
local automorphisms (i.e. those fixing the given point) of finite
type hypersurfaces in $\mathbb C^2$ are linear in some normal
coordinates. The symmetries are then immediately visible from the
defining equation (in the nondegenerate case it follows from a
result of Kruzhlin and Loboda \cite{KL}, in the degenerate case
from \cite{Ko2}.
 See also Section 8 below).

On the other hand, for automorphisms not fixing the point this is
no longer true.
 The
simplest example is given by   rigid hypersurfaces (admitting a
transversal infinitesimal CR automorphism). The normal forms
mentioned above do not respect this symmetry.

For local analysis on  domains which admit symmetries not fixing
the boundary point,
 it is desirable to have  a normalization  which reflects
the symmetries.   This problem was first considered by N. K.
Stanton (\cite{S}), who considered the local equivalence problem
for rigid hypersurfaces of finite type in $\mathbb C^2$ and
constructed a rigid normal form. The results of \cite{S} describe
 all transformations
preserving the rigid normal form and give a complete
classification of rigid hypersurfaces, provided that the model is
not a tube.

In this  paper we consider real analytic Levi degenerate
hypersurfaces of finite type with a tube model. In view of
Stanton's results, this is the only  case of further interest. On
the one hand,
 any hypersurface
which admits a transversal infinitesimal CR vector field is
necessarily rigid. On the other hand, if it admits a
nontransversal one, its model has to be a tube.

The case of real analytic tubes is interesting also in connection
with  the work of G. Francsics  and N. Hanges. In \cite{FH}  they
analyze boundary behaviour of the  Bergman kernel for Levi degenerate
tubes. In relation to this work, Nicholas Hanges formulated the
problem of describing all biholomorphic maps between tubes (\cite{H}).

In recent years, the local equivalence problem on Levi degenerate
hypersurfaces  has been intensively studied (see e.g. \cite{ELZ},
\cite{E}, \cite{Ko1}). In particular, we mention the result of Kim
and Zaitsev (\cite{KZ}), which shows that the second, intrinsic
approach of Cartan, Chern and Tanaka is in general not available.
There has been  substantial progress in understanding the problem
also for CR manifolds of higher codimension (e.g. \cite{Sl},
\cite{GM}, \cite{Sp})

  After introducing notation, we define
 in Section 3  a complete tubular normal form for a general hypersurface
 with a tube model.
It gives the main tool for analyzing biholomorphisms of symmetric
hypersurfaces. The construction is analogous to that of \cite{Ko1}
and is given on the level of fomal power series. The fact that
such a construction can be used for classification problems
 relies
on the essential result of Baouendi, Ebenfelt and Rothschild on
convergence of formal equivalences (\cite{BER1}).
  Rigid
hypersurfaces are considered in Section 4, where a
rigid normal form is obtained. Then we analyze biholomorphisms
preserving this normalization.

In Section 5 we consider tubes. The  symmetry preserving
biholomorphisms are shown to be linear, described by three real
parameters.
 In particular, we obtain an answer to the
question of Hanges. Further we prove that the complete normal form
of Section 3 is convergent for all tubes, thus providing a
complete, convergent and symmetry preserving normal form.

 For non-tubular rigid hypersurfaces we show in Section 6
that only a one parameter family of dilations preserves the rigid
normal form. We apply this result to show that  Stanton's normal
form (\cite{S}) has the same property, and provides therefore a
complete, convergent and symmetry preserving normal form for the
class of rigid hypersurfaces.

Section 7 we consider hypersurfaces which admit
nontransversal infinitesimal CR automorphisms and define a
complete normalization for this class of hypersurfaces.

In Section 8 we give a proof of linearity of local automorphisms
in the normal coordinates of \cite{Ko1}. As a consequence, we
obtain a refinement of the classification result of \cite{Ko2}. It
applies to hypersurfaces with finite local automorphism group, and
allows to determine immediately the size of this cyclic group from
the defining equation in normal coordinates.

As already mentioned, the normal coordinates are a priori only
formal. On the other hand, by the classification result, the local
automorphism group is noncompact if and only if the hypersurface
is a model, when the local automorphisms are  linear already in
the canonical coordinates. In the remaining cases, when the group
is compact, it follows from Bochner's theorem (\cite{B},
\cite{MZ}) that there exist genuine (convergent) coordinates in
which the local automorphisms are linear.

Part of this work was done when the author was visiting the
follow-up program Complex Analysis,
Operator Theory, and Applications to Mathematical Physics in
ESI. He would like to thank Friedrich Haslinger for the invitation and
hospitality, and for the support recieved from ESI.

\section{Preliminaries}

Let $M\subseteq \Bbb C^2$ be  a real analytic hypersurface and
$p\in M$ be a point  of finite type $k$ in the sense of J. J. Kohn
(\cite{K}).

 We will describe $M$ in a neighbourhood of  $p$ using local
holomorphic coordinates $(z,w)$ centered at $p$, where  $z = x +
iy,\  w = u + iv$. The hyperplane $\{ v=0 \}$ is assumed to be
tangent to $M$ at $p$. $M$ is then described near $p$ as the graph
of a uniquely determined real valued function
$$ v = F(z, \bar z, u).$$
 Recall that $p \in
M$ is a point of finite type
if and only if there exist local holomorphic coordinates such that
$M$ is given by
\begin{equation}v = \sum_{j=1}^{k-1}
a_j z^j\bar z^{k-j} + o(\ab k, u ),\ \label{2.9}\end{equation}
where the leading term  is a  nonzero real valued homogeneous
polynomial of degree $k$, with $a_j \in \Bbb C$ and $a_j =
\overline{a_{k-j}}.$

The model hypersurface to $M$ at $p$ is defined using the leading
homogeneous term,

\begin{equation} M_H = \{(z,w) \in \cdva\ | \
 v  = \sum_{j=1}^{k-1}
a_j z^j\bar z^{k-j} \}. \label{2.10}\end{equation}
 When the leading term is equal to  $ \ab k$, we will write
\begin{equation}S_k = \{(z,w) \in \Bbb C^2 \ \vert \   v = \ab k
\}. \end{equation}

Two basic integer valued invariants used in the normal form
construction in \cite{Ko1} will be needed. The first one, denoted
by $e$, is the essential type of the model hypersurface. It can be
described as the lowest index in (\ref{2.9}) for which $a_{e} \neq
0$.

When $e < \frac{k}2$, the second invariant is defined as follows.
Let $e=m_0 <m_1 <\dots <m_s <\frac k2$ be
 the indices in (\ref{2.9}) for which $a_{m_i}\neq 0$. The invariant,
 denoted
 by $L$,
 is the greatest common divisor of the numbers
$\ k-2m_0, k-2m_1,  \dots, k-2m_s$.

It was proved in \cite{Ko1} that for $e < \frac{k}2$ the local
automorphism group $Aut(M_H,p)$  of $M_H$ consists of linear
transformations
$$z^* =   \delta \exp {i\theta} z,\ \  \ \ \   w^* = \delta^k w,$$
                           where
$\exp {i\te}$ is an $L$-th root  of unity and $\delta > 0$ for $k$
even or $\delta \in \Bbb R \setminus \{ 0 \}$ for $k$ odd.

 For $e = \frac{k}2 $ the
local automorphism  group of $S_k$ has dimension three. Its
elements are of the form $z^* = \tf(z,w),\ \ w^* =  \tg(z,w)$,
where
\begin{equation}
\tf(\delta, \mu, \theta; z,w) =
   \frac{ \delta \exp {i\theta} z}{(1 + \mu w)^{\frac1e}}, \ \ \
\ \tg(\delta, \mu, \theta; z,w)=
  \frac{ \delta^k w}{1 + \mu w},\label{4.11}\end{equation}
for $\delta > 0,$ and $\te, \mu  \in \mathbb R$.

\section{a tubular normal form}

We will assume that the  model at $p\in M$ is a tube. By
appropriate scaling and adding a  harmonic term we may assume that
the leading term is equal to $(\frac{ z + \bar z }2)^k$. In
particular, $e$ is equal to one. The model hypersurface is now
\begin{equation}T_k = \{(z,w) \in \Bbb C^2 \ \vert \   v = x^k
\}. \end{equation} $Aut(T_k,0)$ is isomorphic to $\mathbb R^*$,
and consists of dilations
$$ z^* = \delta z, \ \ \ \ \ w^* = \delta^k w,$$
where $\de \in \mathbb R^*$.

 A standard weight assignment will be used. The
variables $z,x,y $ are given weight
 one and  $w$ and $u$ weight $k$.
The defining equation has form
\begin{equation}
v =  F (x,y,u),\label{xk}
\end{equation}
 where
\begin{equation}
 F (x,y,u) = x^k + \sum_{j+l+km\geq k+1} A_{j,l,m}\; x^jy^lu^m.
\label {F}\end{equation} Hence $F(x,y,u) - \;x^k\;$  contains
precisely terms of weight greater or equal to $k+1$.
 Consider the partial Taylor expansion of $F$ in $x,y$
and write
$$F (x,y,u) = x^k + \sum_{j,l} X_{jl}(u) x^j y^l, $$
where $X_{jl}(u) = \sum_{m} A_{j,l,m}\; x^jy^lu^m $.

 We decompose $F$ into parts containing terms of equal
weight:
$$F = x^k +  \sum_{\nu=k+1}^{\infty} F_{\nu},$$
and subject $F$ to a  transformation of the form
\begin{equation}\begin{aligned}
  z^* &= z + f(z,w)\\
  w^* &= w + g(z,w), \\
\end{aligned} \label{fg} \end{equation}
where
\begin{equation}\begin{aligned}
f(z,w) = \sum_{wt. > 1} & f_{jm} z^j w^m\\
g(z,w) =  \sum_{wt. > k} & g_{jm} z^j w^m. \\
\end{aligned} \label{ffg} \end{equation}

 Such transformations preserve form (\ref{xk}), (\ref{F}).
Conversely, it's easy to verify that any transformation preserving
form (\ref{xk}), (\ref{F}),  can be written uniquely as a
composition of an element of $Aut(T_k,0)$ and a transformation of
this form (this factorization will be used repeatedly in the
sequel).

 Again
we decompose the power series
 into parts of the same weight
$$f=\sum_{\nu = 2}^{\infty} f_{\nu} \ \ \ \ \
\text{and} \ \ \ \ \ g=\sum_{\nu =k+1}^{\infty} g_{\nu},$$ and
denote such a transformation by $(f,g)$. Let $v^* = F^*(x^*, y^*,
u^*)$ be the new defining equation, where
 $$F^*(x^*, y^*, u^*)=(x^*)^k + \sum_{j+ l+mk \geq k+1} A^*_{j,l,m}
 (x^*)^j(y^*)^l(u^*)^m.
$$
We will  have to consider also formal hypersurfaces and formal
transformations. From now on
 we  allow both $F, F^*$ and $f,g$ to be formal power series.
The power series formulae  are then interpreted in this  sense.

 Substituting (\ref{fg}) into $v^* = F^*(x^*, y^*, u^*)$
and restricting the variables to $M$,  we get the transformation
formula
\begin{equation}
\begin{aligned}
F^*(x + Re\; f(x+iy, u + i F(x,y,u)),y + Im&\; f(x+iy, \\
 u + i F(x,y,u)),  u+ Re\; g( x+iy, u + i F(x,&y,u)) = \\
 Im \; f(x+iy,
u + i F(x&,y,u))
 + F(x,y,u)
.\label{cov}\end{aligned}\end{equation} In principle, by
multiplying out one can obtain equations for coefficients of
$F^*$, expressed in terms of $F, f, g$. The group of formal
transformations (\ref{fg}) acts on formal power series
 (\ref{F}) via this transformation formula.
\\[2mm]
\defi{We say that $F$ is in t-normal form if
\begin{equation}
X_{0,j} = X_{1,j} = X_{k-1, j} = X_{k, j} = 0,
\label{aaa}\end{equation} for all $ j = 0, 1, \dots,$ and
\begin{equation}
X_{2k-1, 0} = X_{2k-1,1} = 0. \label{ab}\end{equation} }
\\[2mm]
 \pro{ There is a unique
 formal
transformation  (\ref{fg}) which takes $M$ into t-normal form.}
\\[2mm]
{\it proof.} Using induction on weight we will prove that the
normal form conditions on $F^*$ determine uniquely all
coefficients of $f$ and $g$ in (\ref{fg}).
 Let us consider terms of weight
$\mu > k$ in (\ref{cov}). We have
$$\begin{aligned} F^*_{\mu}(x, y, u) +
k Re \; x^{k-1}  f_{\mu-k+1}(x+iy,u + &ix^k)= \\
 = F_{\mu}(x,y, u) + Im\ &g_{\mu}(x+iy,u+ix^k) + \dots
 \end{aligned}$$
where dots denote terms depending on $f_{\nu-k+1}, g_{\nu},
F_{\nu}, F^*_{\nu}$ for  $\nu < \mu$.
 We denote
\begin{equation}L_{\mu}(f,g) =  Re \{ig_{\mu}(x+iy,u+i x^k) + k x^{k-1} f_{\mu - k
+1}(x+iy,u + i x^k)\}, \label{lfg}
\end{equation}
 an analog of the Chern-Moser operator.
For individual monomials in (\ref{lfg}) we have
\begin{equation}\begin{aligned}
k x^{k-1} Re\; \{f_{jm} (x+iy)^j(u+&ix^k)^m\} =
 k x^{k-1}Re\; \{ f_{jm}
(x+iy)^j\}
 u^m \\- m k x^{2k-1} &Im\; \{f_{jm}(x+iy)^j \}
 u^{m-1} + O(x^{2k})\label{ff}\end{aligned}
\end{equation}
and
\begin{equation}\begin{aligned}
Re\; \{i g_{jm} (x+iy)^j(u+i&x^k)^m\} = - Im\; \{g_{jm} (x+iy)^j\}
 u^m \\ - m x^{k} &Re\; \{g_{jm}(x+iy)^j\}
 u^{m-1} + O(x^{k+1}).\label{gg}\end{aligned}
\end{equation}

 In this expansion we will collect coefficients of $x^jy^lu^m$
in (\ref{cov}). Denote $B^*_{j,l,m} = A_{j,l,m} - A^*_{j,l,m}$.
  First we consider $j=0$ and $j=1$.
 Since $k>2$,
 all terms in (\ref{ff}), (\ref{gg}) are multiples of $x^2$,
except for the first  term in (\ref{gg}). That gives for $l = 1,2,
\dots$
\begin{equation}\begin{aligned}
B^*_{0,0,m} & =  - Im \; g_{0,m} + \dots \\
B^*_{0,l,m} & =  - Im \; (i^l g_{l,m}) + \dots\\
B^*_{1,l-1,m} & =- l\; Im\; (i^{l-1}  g_{l,m}) + \dots,  \\
\end{aligned}
\label{11}\end{equation}  where dots denote terms depending on
$f_{\nu-k+1}, g_{\nu}, F_{\nu}, F^*_{\nu}$ for  $\nu < \mu$, which
have been already determined.  Hence the condition that
$A^*_{0,0,m}=0$ determines $Im\; g_{0,m}$ and $A^*_{0,l,m}=0$
determines $g_{l,m}$ for $l=1,2,\dots$. Further we consider
$j=k-1$ and $j=k$. For $l \geq 2$,   we get contribution from the
first term in (\ref{ff}) and the  two terms in (\ref{gg}). This
gives
\begin{equation}
\begin{aligned}B^*_{k-1,l,m} & =    k Re
\; (i^{l} f_{l,m}) -
 C_1 Im\;  ( i^{l} g_{k+l-1,m}   ) + \dots \\
B^*_{k,l-1,m} & =  k l\;Re\; (  i^{l-1} f_{l,m}) - C_2 Im \;
(i^{l-1} g_{k+l-1,m}) -  \\  &(m+1) Re\; (i^{l-1} g_{l-1,m+1})
\dots , \\
\end{aligned}
\label{12}\end{equation}
 where
$C_1 = \binom{k+l-1}{l}$ and  $C_2 = \binom{k+l-1}{l-1}$.
 That determines $f_{l,m}$ for
$l = 2,3, \dots$. Next we consider $(j,l) = (k-1,0)$, with
contribution from
 the first terms
in (\ref{ff}) and (\ref{gg}).
 For $(j,l) = (2k-1,0)$
 all terms
contribute, and we obtain
\begin{equation}
\begin{aligned}
B^*_{k-1,0,m} & =  k Re \ f_{0,m} -
Im\; g_{k-1,m}+ \dots \\
B^*_{2k-1, 0,m-1} & = k Re f_{k, m-1} - km Im f_{0,m} \\
 & - Im \;
g_{2k-1,m-1} - m Re \; g_{k-1,m} + \dots.
\end{aligned}\label{13}\end{equation}
That  determines $f_{0,m}$, since all other entries have been
already determined.  For $(j,l) = (k-1,1)$
 the  first terms
in (\ref{ff}) and (\ref{gg})
 contribute, for $(j,l) = (2k-1,1)$ all contribute.
That gives
\begin{equation}
\begin{aligned}
B^*_{k-1,1,m} & =   - k Im\; f_{1,m}
- k Re\;  g_{k,m} + \dots  \\
B^*_{2k-1,1,m-1} & =  - k (k+1)Im \; f_{k+1,m-1}
  -
  kmRe\;  f_{1,m} \\ & - 2k Re \;
g_{2k,m-1} + k m Im \; g_{k,m},
 \end{aligned}\label{14}\end{equation}
which determines  $f_{1m}$. For $(j,l) =  (k,0)$ we get
contribution from the first term in (\ref{ff}) and the two terms
in (\ref{gg}), which gives
\begin{equation}
B^*_{k,0,m}  =   k Re \; f_{1,m} - Im\; g_{k,m}
 -  (m+1) Re\; g_{0,m+1} + \dots. \label{15}\end{equation}
That determines $Re \; g_{0,m+1}$. It is immediate to verify that
the first appearance of each equation agrees with the
normalization conditions on $f$ and $g$.
\\[2mm]
{\bf Remark 3.1.} Note that the same result is obtained if the
conditions $X_{2k-1,0}= X_{2k-1,1} = 0 $  are replaced by
$X_{2k-1,0} =A$ and  $X_{2k-1,1}=B$, for any fixed real numbers
$A,B$. This remark will be used in the proof of Proposition 6.1.
\\[2mm]
From the factorization of a general map preserving form
(\ref{xk}), (\ref{F}), we obtain the following corollary.
\\[2mm]
{\bf Corollary 3.1} {\it \; The only transformations which
preserve the t-normal form are the elements of $Aut(T_k,0)$.}

\section{A rigid normal form}
In this section we consider rigid hypersurfaces with tube models
and define a rigid t-normal form.

Consider a rigid hypersurface with a tube model, given by

\begin{equation}
v = x^k + F (x,y),
\end{equation}
where
\begin{equation}
 F (x,y) =\sum_{j+l \geq k+1} A_{j,l} x^jy^l,
\end{equation}
and $k>2$.
\\[2mm]
\defi{ We say that $F$ is in rigid t-normal form if

\begin{equation}
A_{0,j} = A_{1,j} = A_{k-1, j} = A_{k, j} = 0
\label{aa}\end{equation} for all $ j = 0, 1, \dots$.}
\\[2mm]
\lema{ There exists a unique transformation of the form
\begin{equation}\begin{aligned}
  z^* &= z + \sum_{i=2}^{\infty} f_{i} z^i\\
  w^* &= w + \sum_{i=k+1}^{\infty}  g_{i} z^i, \\
\label{4.1}
\end{aligned} \end{equation}
which takes $F$ into  rigid t-normal form.}
\\[2mm]
{\it proof:}
 We will again determine by induction the coefficients $f_i$, $g_i$
 in such a way that (\ref{aa}) is satisfied.
 Let
\begin{equation}
v^* = (x^*)^k + F^* (x^*,y^*)
\end{equation}
in the new coordinates,  where
\begin{equation}
 F^* (x^*,y^*) =\sum_{j+l \geq k+1} A^*_{j,l} (x^*)^j(y^*)^l.
\end{equation}
The transformation formula takes form
\begin{equation}
\begin{aligned}
(  x + Re f (x+iy))^k& + F^*(x + Re\; f(x+iy), y + Im \; f(x+iy) )
\\ = Im \; &g(x+ iy)\;  +\;  x^k + F(x,y).
\end{aligned}
\end{equation}
Terms of degree $m > k$ in this equation depend linearly on terms
of degree $m-k+1$ in $f$ and degree $m$ in $g$, and nonlinearly on
terms of lower degree in $f$. For the terms specified in
(\ref{aa}) we have  the following equations
\begin{equation}\begin{aligned}
B^*_{0,l} & = - Im\; (i^l g_l)  + \dots\\
B^*_{1,l-1} & = - l\; Im\; (i^{l-1} g_l) + \dots,
 \end{aligned}\end{equation}
 and
\begin{equation}\begin{aligned}
B^*_{k-1,l} & = k Re\; (i^{l}  f_{l}) -
 C_1 Im\; ( i^{l} g_{l+k-1}) + \dots\\
B^*_{k,l-1} & = k l  Re \;(i^{l-1} f_{l})
  - C_2 Im\; (i^{l-1}
g_{l+k-1}) + \dots,
 \end{aligned}\end{equation}
where $B^*_{j,l} = A_{j,l} - A^*_{j,l}$, $C_1 = \binom{k+l-1}{
l},$ $C_2 = \binom{k+l-1}{l-1}$, and dots denote already
determined numbers. The first two equations determine $g_l$, $ l =
k+1, k+2,  \dots$,  the second two determine $f_l$, $l=2, 3,
\dots$.
 \qed

 Using the
complete t-normal form we will analyze biholomorphic
transformations which preserve  the rigid t-normal form.

\section{Biholomorphic equivalence of tubes}
Since a tube hypersurface satisfies automatically all t-normal
form conditions except for $A_{2k-1,0,0} = 0$, this normalization
can be used effectively to classify tubes.
\\[2mm]
 \pro{ Let $M_1$ and $M_2$ be two tubular hypersurfaces of finite
  type $k>2$ at the origin, given
by $v = F(x)$ and $v = G(x)$, respectively. If $\Psi$ is a local
biholomorphism preserving the origin, which maps $M_1$ to $M_2$,
then it has form
$$\begin{aligned}
  z^* = az + ibw, \ \ \ \ \ \ \
  w^* = cw  \\
\end{aligned} $$
 for $a,c \in \mathbb R^*$ and $b \in \mathbb R$. In this case,  $F$ and $G$ satisfy
$$ G(ax   - bF(x)) = c F(x).$$
}
\\[2mm]
{\it  proof.} By assumption, $F(x) = c_1 x^k + o(x^k)$ and $G(x) =
c_2 x^k + o(x^k)$ for some nonzero real constants $c_1, c_2$.
First we make those coefficients equal to one, using the dilations
$w^* = c^{-1}_j w$, $j = 1,2$.
 In the second step, we put $M_1$
into t-normal form by a transformation
$$ z^* = z + ihw, \ \ \ \ \ w^* = w.$$
In (\ref{cov}) we have $Re\;(ih(u+iF)) = -hF$,  so the
transformation equation becomes
\begin{equation} F^*(x - hF(x)) = F(x). \end{equation}
 Equating
coefficients of $x^{2k-1}$ we obtain $ A^*_{2k-1,0,0} - kh =
A_{2k-1,0,0} $. Hence for $ h = -\frac1k A_{2k-1,0,0}$ the
t-normal form is obtained. Next we perform the same normalization
on
 $M_2$,
 and  consider
the two resulting  hypersurfaces in t-normal form. By Corollary
3.1., a biholomorphic equivalence between them is an element of
$Aut(T_k,0)$. By composing the five linear mappings, we obtain the
claimed form of the biholomorphism.
 The transformation formula gives
$G(ax - bF(x)) = c F(x).$
\\[2mm]
By the same reasoning, we obtain the following corollary.
\\[2mm]
 {\bf Corollary 5.1.} {\it \ For any tube hypersurface the
complete t-normal form is convergent  and preserves the tubular
symmetries.}

\section{Nontubular hypersurfaces}

The t-normal form preserves tubular symmetries, but does not in
general preserve rigidity.
 In this
section we determine biholomorphisms which preserve the rigid
t-normal form for hypersurfaces other than tubes. Then we obtain
the same result for  Stanton's normal form.
\\[2mm]
\pro{\it Let $M$ be a rigid hypersurface, which is not equivalent
to a tube. Then the only transformations which preserves the rigid
t-normal form are the
 dilations
\begin{equation}
  z^* = \delta z \ \ \ \ \
  w^* = \delta^k w, \label{ac}
\end{equation}
where $\delta \in \mathbb R^*$.}
\\[2mm]

 {\it proof:}
 Let $M_1$, $M_2$   be two hypersurfaces
in rigid t-normal form, given by $v = F(x,y)$ and
  $ v^* =
F^*(x^*,y^*)$, respectively. $F$ satisfies all t-normal form
conditions, except possibly for $A_{2k-1,0,0}= 0$ and
$A_{2k-1,1,0} = 0$, and the same holds for $F^*$.  Let $\Psi$ be
  a
local biholomorphism preserving the origin, which maps $M_1$ to
$M_2$. Since $\Psi$ preserves form (\ref{xk}), (\ref{F}), we write
it again as the composition of an element of $Aut(T_k,0)$ and a
transformation of the form
$$\begin{aligned}
  z^* = z + \sum_{wt. > 1}  f_{jm} z^j w^m,\ \ \ \ \
  w^* = w + \sum_{wt. > k}  g_{jm} z^j w^m.
\end{aligned} $$
We consider this transformation and show that all the coefficients
of $f$ and $g$ have to vanish.
  The transformation formula now takes form

\begin{equation}\begin{aligned}
F^*(x + Re f(x+iy,& u + i F(x,y)), y + Im f (x+iy, u + i F(x,y))
\\ = Im \; g(x+iy, &u + i F(x,y)) + F(x,y).\label{xiy}
\end{aligned}\end{equation}
 Let
$$F(x,y) = H(x) + Q(x,y) + o_{wt}(p), $$
where
$$Q(x,y) = \sum_{j=1}^{p} d_{j} y^j x^{p-j}. $$
Hence $Q$ is the first homogeneous level containing $y$. Let $j_0$
be the first index for which $d_{j_0}$ is nonzero. Similarly, let
$$F^*(x^*, y^*) = H^*(x^*) + Q^*(x^*,y^*) + o_{wt}(p^*), $$
with $Q^*(x^*,y^*) = \sum_{j=1}^{p^*} d^*_{j} (y^*)^j
(x^*)^{p^*-j}. $ Without any loss of generality, we assume that $p
\leq p^*$.

Consider (\ref{cov}), (\ref{xiy}) and the resulting equations for
coefficients of each monomial $z^j \bar z^l u^m$. We will  use the
notation introduced in the proof of Proposition 3.1.

 First, since
$F$ and $F^*$ satisfy all t-normal form conditions for weight less
than $2k-1$, it follows from Proposition 3.1. that all terms in
$f$ of weight less than $ k$ and  in $g$ of weight less than
$2k-1$ are zero. In particular, $f_{20}, ..., f_{k-1,0}$ and
$g_{11}, ..., g_{k-2,1}$ all vanish. Next we consider equations
for terms of weight $2k-1$ and $2k$. By (\ref{11}), (\ref{12})  we
obtain $f_{k,0} = f_{k+1,0} = 0$ and $g_{k-1,1} = g_{k,1} =
g_{2k-1,0} = g_{2k,0} = 0$. By rigidity and (\ref{13}),
(\ref{14}), we get  $Re\; f_{0,1} = 0$ and  $Im\; f_{1,1} = 0$.
Hence the only possibly nonzero coefficients in $f$ of weight less
than or equal to $k+1$ and in $g$ of weight less than or equal to
$2k$ are $Im\; f_{0,1}$ and $Re\; f_{1,1}$.
 Further,
consider equations for terms in the range of weights $2k+1, \dots,
p+k$. Terms in $F$ and $F^*$ of weight greater than  $p$ enter
these equations only through $f_{20}, ..., f_{k,0}$ and $g_{11},
..., g_{k-1,1}, $ which we already know to be zero. So the
equation for the coefficient of $ x^{p - j_0} y^{j_0-1} u$ comes
only from $Q$, namely
\begin{equation}B^*_{ {p - j_0 },{j_0-1},1} = j_0 Im
f_{01}.\label{tuk}
\end{equation}
It follows from rigidity that  $Im f_{01} = 0$. From (\ref{14}),
the equation for the coefficient of $x^{j_0} y^{p - j_0 } u$ is

\begin{equation}B^*_{{j_0}, {p - j_0 },1} = -k (p-j_0) Re\;
f_{11},\label{tuk2}
\end{equation}
so $Re\; f_{11}=0$. From (\ref{13}) and (\ref{14}) we have

\begin{equation}
\begin{aligned}
B^*_{2k-1, 0,0} & = k Re f_{k, 0} - km Im f_{0,1} \\
 & - Im \;
g_{2k-1,0} - m Re \; g_{k-1,1} + \dots.
\end{aligned}\end{equation}

and

\begin{equation}
\begin{aligned}
B^*_{2k-1,1,0} & =  - k (k+1)Im \; f_{k+1,0}
  -
  kmRe\;  f_{1,1} \\ & - k Re \;
g_{2k,0} + k m Im \; g_{k,1},
 \end{aligned}\end{equation}
All terms on the  right hand sides are zero, hence $A_{2k-1, 0,0}
= A^*_{2k-1, 0,0}$ and $A_{2k-1, 1,0} = A^*_{2k-1, 1,0}$. By
Proposition 3.1. and Remark 3.1., all coefficients of $f$ and $g$
vanish,  and $\Psi $ is an element of $Aut(T_k,0)$.

 \qed

We use the above result to show that  Stanton's normal form of a
nontubular hypersurface with a tube model is preserved only by
dilations. Hence this normal form  provides a convergent, symmetry
preserving complete normalization for the class of rigid
hypersurfaces.

Recall that Stanton's normal form for a hypersurface with a tube
model uses the complex Taylor expansion of $F$,
$$v = x^k + \sum_{j,l}A_{jl}z^j \bar z^l.$$
The normal form conditions are
$$ A_{0,l}=A_{1,l}=0,$$
for all $l = 1,2,\dots$
\\[2mm]
{\bf Proposition 6.2.}{\it \ Let $M_1, M_2$ be two nontubular
rigid hypersurfaces in Stanton's normal form and let $\Phi_1$
 be a local biholomorphic map preserving the origin, which maps $M_1$ to $M_2$.
Then $\Phi_1$ is a dilation of the form (\ref{ac}).}
\\[2mm]
{\it proof:}  By Lemma 4.1., there is a formal  transformation
$\Psi_1$ of the form (\ref{4.1}), which takes $M_1$ into rigid
t-normal form, and a formal  transformation $\Psi_2$ of the same
form which takes $M_2$ into rigid t-normal form. We denoted by
 $M_1^T$ and
$M_2^T$ the corresponding  hypersurfaces (apriori only formal).
Then $\Phi_2 = \Psi_1^{-1} \circ \Phi_1 \circ \Psi_2$ is a formal
equivalence of $M_1^T$ and $M_2^T$. By Proposition 6.1., $\Psi_2$
is a dilation. Since $\Phi_1 = \Psi_1 \circ \Phi_2 \circ
\Psi_2^{-1}$, it follows that $\Phi_1$ has form (\ref{4.1}).  By
the results of \cite{S}, the only transformations of this form
which preserve Stanton's normal form are dilations of the form
(\ref{ac}). \qed
\\[2mm]
By the same argument we obtain the following corollary.
\\[2mm]
 {\bf Corollary 6.1. }{\it \
For nontubular rigid hypersurfaces the only transformations
preserving rigidity are transformations of the form (\ref{4.1})}.

\section{nontransversal symmetries}

We now consider  a hypersurface $M$ which admits a nontransversal
infinitesimal CR automorphism. By straightening the corresponding
vector field, we obtain local holomorphic coordinates in which the
defining equation has  form

\begin{equation} v = x^k + G(x,u), \label{vg}\end{equation}
where $G(x,u)$  is $o_{wt}(k)$, and
$$G(x,u) = \sum_{j=0}^{\infty} X_j(u) x^j.$$
\\[2mm]
\defi{
$M$ is in normal form if the defining equation  has form
(\ref{vg}), and satifies
\begin{equation}
 X_{0} =X_{k-1}=
X_{k}=  X_{2k-1} = 0. \label{a}\end{equation} }
\\[2mm]
\pro{There exists a transformation of the form
\begin{equation}
z^* = z + \psi(w), \ \ \ \ \  w^* = w + \phi(w) \label{ps}
\end{equation}
which takes $M$ into normal form.}
\\[2mm]
{\it proof:} First we  show that the  transformations of the form
(\ref{ps}) preserves form (\ref{vg}). We have \begin{equation}
\begin{aligned}
 G^*(x + Re\; &\psi(u + iG(x,u)), y + Im \;
\psi(u + iG(x,u), u) = \\& G(x,u) +   Im \; \phi(u + i G(x,u)).
\end{aligned}
\end{equation}
Since the right hand side is independent of $y$, it follows
immediately that $G^*$ is  independent of $y^*$.
 By  (\ref{11}) - (\ref{14}),  the
equations for $B^*_{0,0,m}$, $B^*_{k-1,0,m}$ $B^*_{k,0,m}$ and
$B^*_{2k-1,0,m}$ involve precisely the coefficients  $f_{j,m}$ and
$g_{j,m}$ with $j=0$.
 Hence, setting all other coefficients equal to zero,
  the induction argument of Proposition
3.1. determines the coefficients $f_{0,m}$ and $g_{0,m}$ so that
(\ref{a}) is satisfied.
\qed
\\[2mm]
 Now we show that the conditions
in Definition 7.1. define a complete normalization.
\\[2mm]
 \pro{
Let $\Psi$  be a transformation which preserves (\ref{vg}) and the
normal form conditions (\ref{a}). Then $\Psi$ is a dilation
$$z^* = \delta  z
, \ \ \ \ \ w^* = \delta^kw.$$}
\\[2mm]
{\it proof:\ } We  decompose again $\Psi$ into an element of
$Aut(T_k,0)$ and a transformation $(f, g)$ of the form (\ref{fg}),
and consider the effect of this transformation.  Assume $(f, g)$
is not the identity and let $\mu$ be the first weight where it
differs from the identity. The defining equation of $M$ satisfies
all t-normal form conditions except possibly for $X_{10} = 0$.
Hence, as in Remark 3.1., the equations for weight $\mu$  have to
contain the one for $x u ^m$, with  some $m
>0$, and $B^*_{1,0,m} $  has to be nonzero.
Consider the equation for  the coefficient of $x^{2k-1} y^2
u^{m-2} $.  From \ref{ff} and \ref{gg} we obtain
 \begin{equation}
\begin{aligned}
 0 &= - k \binom{k+2}{2} Re\;  f_{k+2,m-2} +
 m k Im \; f_{2,m-1}\\ & +  \binom{2k+1}{2} Im\;g_{2k+1,m-2}
+ m \binom{k+1}{2} Re \; g_{k+1,m-1}.
\end{aligned}\end{equation}
 By (\ref{11}) and (\ref{12}), we have $Re
g_{k+1,m-1} = Im g_{2k+1,m-2} = Re \; f_{k+2,m-2} = 0$. Hence also
$Im \; f_{2,m-1} = 0$. Now consider the second equation in
(\ref{12}), with $l=2$. It gives $Im \;g_{1,m} = 0$. So by
(\ref{11}), the equation for $x u ^m$ is $B^*_{1,0,m}=0, $ which
is a contradiction. Thus we proved that the only transformations
preserving the normalization conditions are the elements of
$Aut(T_k,0)$.
\qed

\section{Linearity of local automorphisms}

In this section we consider a general hypersurface of finite type
and prove linearity of local automorphisms in normal coordinates.

We will use also Taylor expansion of $F$ in terms of $z, \bar z,
u$,
$$F(z, \bar z, u) = \sum_{j,l}Z_{jl}(u)z^j\bar z^l,$$
where
 $$Z_{jl}(u) = \sum_{m} a_{jlm} u^m,$$
 and consider the complete normal forms obtained in \cite{Ko1}.
The following result in the case $e=\frac k2$ was obtained in
\cite{Ko2}.
\\[2mm]
 \pro{\it Let $M$ be a Levi degenerate hypersurface of finite type,
  not equivalent to $v = |z|^k$. Then
 all local automorphisms expressed in normal
coordinates are linear.}
\\[2mm]
 {\it proof:}
First we consider the case $e<\frac k2$. If $M_H$ is different
from a tube,  the normal form conditions are
 $$\begin{array} {rl}
Z_{j0} & = 0, \ \ \ \ \ j=1,2,\dots,  \\
Z_{k-e+j,e} & = 0, \ \ \ \ \ j= 0,1,\dots, \\
Z_{2k-2e, 2e} & = 0, \\
(Z_{k-1}, P_z) & = 0,  \label{3.11} \end{array}.$$ where
$$(Z_{k-1}, P_z) = \sum_{j=1}^{k-2}Z_{j,k-1-j}  (j+1)\bar a_{j+1}.$$

 The symmetry
group of the model acts on normal forms. We will prove that each
element of $Aut(M_H,p)$ preserves the normal form, hence its
action on normal forms is direct (no renormalization is needed).
Since every element of $Aut(M_H,p)$ is linear, its  application
clearly preserves the first three conditions. In order to see that
the last condition is also preserved, we write $P$ as
$$ P( z, \bar z ) = \sum_{j=1}^c a_{jL} z^{jL}\bar z^{(c-j)L},$$
where $c = \frac{k}{L}$, and
$$Z_{k-1}(u)  =  \sum_j  \beta_j(u) z^{k-1}\bar z^{k-1-j}.$$
We have
$$ P_z(z, \bar z) =\sum_{j=1}^c jL a_{jL} z^{jL-1}\bar
z^{(c-j)L}, $$
 and consider the action of a transformation $z^* = \al z$,
where $\al^L = 1.$ We obtain

$$(Z^*_{k-1}, P_z) = \sum_{j=1}^{c}
jL a_{jl}  \beta_{jL-1}(u) \al^{jL-1} \bar \al^{(c-j)L}
  =
  $$
  $$=\al^{-1} \sum_{j=1}^{c}
jL a_{jl}  \beta_{jL-1}(u)  \bar \al^{(c-j)L}  = (Z_{k-1}, P_z) =
0.$$

Now let the model be a tube. In this case $L=2$ or $L=1$,
depending on the parity of $k$. The normal form from \cite{Ko1}
are

$$\begin{aligned} Z_{j0} &= 0, \ \ \ \ \ j=1,2,\dots,  \\
Z_{k-1+j,l} &= 0, \ \ \ \ \ j= 0,1,\dots, \end{aligned}$$ and
$$
Z_{2k-2, 2} = Re \; Z_{k-2,1} = Re \; Z_{k, k-1}  = 0.
 \label{4.5} $$
Clearly this normalization is preserved when an element of
$Aut(T_k,0)$ is applied.

 Next consider the third case, $e=\frac k2$.
For completeness we repeat the argument here.
 Let us consider normal coordinates for $M$, i.e. F satisfies
$$
Z_{j0} = Z_{e,e+j} = 0
$$
for all $j = 0,1,2, \dots$,  and
$$
Z_{2e, 2e} = Z_{3e, 3e} = Z_{2e, 2e-1} = 0.
$$
We separate the first two  leading terms in the Taylor expansion
of $F$,
\begin{equation} \fz = \ab k + \qz + o_{wt}(p), \end{equation}
where $\tilde P$ is a nonzero weighted homogeneous real valued
polynomial of weight $ p > k $
\begin{equation} \qz = \sum_{j + l  + k m = p} a_{j l m}z^{j} \bar
z^{l} u^{m}, \end{equation} and $o_{wt}(p)$ denote terms which are
of weight greater then $p$. We define the index $(j_0, l_0, m_0)$
to be the smallest one in inverse lexicographic ordering (the last
components are compared first, then the second ones) for which
$a_{j_0 l_0 m_0} \neq 0$.

 Let $(f,g)$ be a local automorphism of $M$,
i.e. a transformation which preserves $F$. Its general form
is
\begin{equation}
f(z,w) =\delta e^{i\theta} z + o_{wt}(1) \ \ \ \ \
 g(z,w) = \delta^k w + o_{wt}(k).
\end{equation}
 The
numbers $\delta, \; \theta \; $ and $ \mu = Re\; g_{ww}$
 are the initial data of the automorphism. We  consider
simultaneously  $M$ with the automorphism $(f,g)$ and the model
$S_k$ with the automorphism $(\tf, \tg)$ having  the same initial
data as $(f,g)$. We will use (\ref{cov}) to compare the
coefficients of $(f,g)$ and $(\tf, \tg)$.

 In two steps we will show that $f$ and $g$ may
be replaced by $\tf$ and $\tg$ when considering terms of weight
less or equal to $p+k$ in (\ref{cov}). More precisely,
\begin{equation}f(z,w) =  \tf(z,w) + o_{wt}(p+1),\ \ \ \ \
 g(z,w) = \tg(z,w) + o_{wt}(p+k).\end{equation}
First, since $\tilde P$ has weight p, all equations obtained from
(\ref{cov}) for coefficients of monomials  up to weight $p-1$ are
the same as those for $S_k$ and $(\tf$, $\tg)$. Hence $f$ is equal
to $\tf$ modulo $ o_{wt}(p-k)$ and $\tg$ equal to $g$ modulo $
o_{wt}(p-1)$. For terms of weight $p$, $\tilde P$ enters
(\ref{cov}) only via the linear part of $(f,g)$, as $\tilde
P(\delta e^{i\theta} z, \delta e^{-i\theta}  \bar z, \delta^k u)$.
 Since $\tilde P$ (and in particular  $a_{j_0,
l_0, m_0}$)  has to be preserved, we obtain immediately that
 $\delta = 1$ and $e^{i(j_0 - l_0 )
\theta} = 1.$
 For terms of weight $p+1, p+2, \dots, p+k$, $\tilde P$
enters (\ref{cov}) only through  the initial data $Re g_{ww}$, and
the coefficients   $f_{20}, \dots f_{k0}$ in $f$ and $g_{11},
\dots  g_{k1}$ in $g$.
 But we
already know these coefficients to be the same as in $(\tf, \tg)$, namely
 zero (if $ k > p-k$
we use  an obvious step by step argument). Since by the result of
\cite{Ko1}  a local automorphism is uniquely determined by its
initial data, it follows that $\tf$ has to agree with $f$ modulo
terms of weight greater than $p+1$ and $\tg$ has to agree with $g$
modulo terms of weight greater than $p+k$.
 This proves the
claim.

Now we consider all terms of weight $k+1, \dots, k+p$ in the
transformation formula (\ref{cov}) . On the right hand side, using
$ g(z,w) = w - \mu w^2 + \dots$
 we
have

\begin{equation} \begin{aligned} &Im \; g(z, u + iF) =
 F - Im\;  \mu (u + i(\ab k + \tilde P +
o_{wt}(p)))^2 + J_1 + \\ +& \; o_{wt}(k+p)  = F + 2\mu u \ab k - 2
\mu u \tilde P + J_1 + o_{wt}(k+p),
\end{aligned}\end{equation}
where  $J_1$  denotes  terms of weight $\leq k+p$ which come only
from $\ab k$, in other words terms which appear in the
corresponding expansion for the model  and $(\tf, \tg)$
  (we will not need this expression  explicitly).
 On the left,
 we get from the leading term
\begin{equation}
\vert f(z, u + i(\ab k + \tilde P + o_{wt}(p)))\vert^k = \vert z -
\frac{\mu}e z( u + i(\ab k + \tilde P + o_{wt}(p)))+
o_{wt}(2k)\vert^k
\end{equation}
which gives
\begin{equation}
 \ab k
- 2\mu Im \; \ab k \tilde P + J_2 + o_{wt}(p+k) =
 \ab k
 + J_2 + o_{wt}(p+k)
,\end{equation}
where again $J_2$  denotes all terms of weight $\leq k+p$ which
come only from $\ab k$.
 From the second term in $F=F^*$ we get
\begin{equation}
\begin{aligned}
 &\tilde P(f,\bar f, Re\ g) = \tilde P(
e^{i\theta} z-e^{i\theta}\frac{\mu}e z(u+ i \ab k)  +
o_{wt}(k+1)), \bar{ \ \ },
\\& u - Re\; \mu (u +
i(\ab k + o_{wt}(k)))^2+ o_{wt}(2k)). \end{aligned}
\end{equation}

By the same argument as we used before for $\tilde P$, since
$f_{20}, \dots f_{k0}$
 and $g_{11}, \dots  g_{k1}$ vanish,
terms of weight greater than $p$ and less or equal to $p+k$  in
$F^*$ influence  (\ref{cov}) only via the
 linear part of $(f,g)$.
Multiplying out and taking into account that terms coming only
from $\ab k$ have to  eliminate each other, we calculate the
coefficients of $z^{j_0} \bar z^{l_0} u^{m_0+1}$ in (\ref{cov}).
We obtain
\begin{equation}
 a_{j_0,
l_0, m_0 +1} - a_{j_0, l_0, m_0}( \frac1e \mu (j_0 + l_0 +l m_0))
 = - 2\mu a_{j_0, l_0, m_0 } +
a_{j_0, l_0, m_0 +1}.
\end{equation}
It will hold if and only if
\begin{equation}
a_{j_0,  l_0,  m_0 } (2\mu - \frac1e \mu ( j_0 + l_0 + l m_0 )) =
0,\end{equation}
 hence
\begin{equation}j_0 + l_0 + l m_0 = k\end{equation}
(recall that $ e = \frac k2$). It follows that either $\mu = 0$,
or $m_0 = 1$ and $j_0 + l_0 = e$. If $m_0 = 1$ we consider the
coefficients of $z^{j_0 + k}\bar z^{l_0 + k}$. From the formulas
above we get
\begin{equation}a_{j_0 + k,l_0 + k,0} + \mu a_{j_0, l_0, 1} = a_{j_0 + k,l_0 +
k,0},\end{equation}  and so $\mu = 0$. Hence there is no
automorphism of $M$ with $\mu \neq 0$, and we proved that every
local automorphism in normal coordinates is linear. \qed

 The
following is a refinement of the classification result obtained in
\cite{Ko2}. It allows to determine immediately the size of the
discrete (cyclic) group of local symmetries.

\pro{\it
 For a given hypersurface
 exactly one of the following possibilities occurs.}
\\
\begin{enumerate} {\it \item $Aut(M,p)$ has real
dimension three. This happens  if and only if $M$ is equivalent to
$v=|z|^k$.
\\
 \item $Aut(M,p)$ is isomorphic to
$\Bbb R^+ \oplus {\Bbb Z}_m.$ This happens if and only if $M$ is a
model hypersurface with $l<\frac k2$
\\
\item $Aut(M,p)$ is isomorphic to $S^1$.  This happens if and only
if the defining equation for $M$  in normal coordinates has form
$$ v = G(\ab 2, u).$$
\item $Aut(M,p)$ is isomorphic to $ {\Bbb Z}_m$ . This happens if
and only if $M$ is not a model hypersurface and in the normal
coordinates $m$ is the largest index such the defining equation of
$M$ can be written in the form $$v = G(\ab 2, z^m, \bar z^m u).$$
}
\end{enumerate}

{\it proof: } Follows immediately from Proposition 8.1, since  all
local automorphisms are linear in normal coordinates and act on
each term separately.
\qed
\\[2mm]


\begin{thebibliography}{99}
\bibitem{BER1}  M.S.Baouendi,
P.Ebenfelt, L.P.Rothschild : \textit{Convergence and finite
determination of formal CR mappings},  J. Amer. Math. Soc. \textbf{ 13}
(2000), p. 697-723.

\bibitem{BER2}  M.S.Baouendi,
P.Ebenfelt, L.P.Rothschild : \textit{Local geometric properties of
real submanifolds in complex space},  Bull.\ Amer.\ Math.\ Soc.\ (N.S.)
\textbf{37} 3
  (2000), p. 309--336.


\bibitem{BB}   E.Barletta, E.Bedford :  \textit{Existence of proper
mappings from domains in $\Bbb C^2$ },  Indiana Univ.  Math. J.
 \textbf{2} (1990), p. 315-338.


\bibitem{BFG}  M.Beals, C.Fefferman, R.Grossman : \textit{Strictly pseudoconvex domains in $\Bbb C^n$},
Bull. Amer. Math. Soc.  \textbf{8} (1983), p. 125-322.


\bibitem{B} V.K.Beloshapka : \textit{On the dimension of the group of
automorphisms of an analytic hypersurface},  Math. USSR, Izv.
\textbf{14} (1980), p. 223-245.

\bibitem{BE} V.K.Beloshapka, V.V.Ezhov : \textit{
Normal forms and model hypersurfaces in $\mathbb C^2$}, preprint.


\bibitem{B} S.Bochner : \textit{Compact groups of differentiable
transformations }, Ann. of Math. \textbf{47} (1945), p. 372-381.

 \bibitem{C1}  E.Cartan :
 \textit{Sur la g\'eom\'etrie pseudo-conforme des hypersurfaces de deux
variables complexes, I }, Ann. Math. Pura Appl.  \textbf{11}
(1932), p. 17-90.


 \bibitem{C2}  E.Cartan : \textit{Sur la g\'eom\'etrie pseudo-conforme
des hypersurfaces de deux variables complexes, II}, Ann.Scoula
Norm. Sup. Pisa   \textbf{1} (1932), p. 333-354.

 \bibitem{CS}  S.-C.Chen and M.-C. Shaw : \textit{Partial differential
equations in Several Complex Variables},  American Mathematical
Society/International Press (2001).

\bibitem {CM} S.S.Chern and J.Moser : \textit{Real hypersurfaces in
complex manifolds},  Acta Math.  \textbf{133} (1974), p. 219-271.

 \bibitem {E}  P.Ebenfelt :  \textit{New invariant tensors in CR
structures and a normal form for real hypersurfaces at a generic
Levi degeneracy}, J. Diff. Geometry  \textbf{50} (1998), p. 207-247.

\bibitem{ELZ}  P.Ebenfelt, B.Lamel, D.Zaitsev :  \textit {Finite jet
determination of local analytic CR automorphisms and their
parametrization by 2-jets in the finite type case},
Geom.\ Funct.\ Anal.\ {\bf  13} 3 (2003), p. 546-573.


\bibitem{EHZ}P.Ebenfelt, X.Huang, D.Zaitsev :
  \textit{The equivalence problem and rigidity for hypersurfaces embedded
  into hyperquadrics}, Amer. J. Math. \textbf{127} (2005), p. 169-191.

\bibitem{E} V.V.Ezhov : {\it
Triviality of scalar linear
type isotropy subgroup by passing to an
alternative canonical form of a hypersurface},
 Complex analysis and applications (Warsaw, 1997),  Ann. Polon. Math.  {\bf 70}
  (1998), p. 85-97.

\bibitem{GM}
T.Garrity, R.Mizner : {\it\ The equivalence problem for
higher-codimensional CR structures} Pacific J. Math. \textbf{177} 2 
(1997),  p.  211--235.

\bibitem{FH}
G.Francsics,  N.Hanges : {\it \ Analytic singularities of the
Bergman
 kernel for tubes}  Duke Math. J. {\bf 108}  (2001), p.  539-580.

\bibitem{H} N.Hanges :  Personal communication.

 \bibitem {Ja}  H.Jacobowitz :
\textit{An introduction to CR structures}, Mathematical Surveys
and Monographs 32, AMS (1990).

\bibitem{KZ} S.Y.Kim, D. Zaitsev : \textit{
Equivalence and embedding problems for CR-structures of any codimension}
Topology {\bf 44}  (2005), p. 557-584.


 \bibitem{K} J.J.Kohn : \textit{Boundary behaviour of
$\bar \partial$ on weakly pseudoconvex manifolds of dimension two},
 J.Diff. Geometry  \textbf{6} (1972), p. 523-542.

\bibitem{Ko1} M.Kol\'a\v r : {Normal forms for hypersurfaces of finite type in
 $ \mathbb C^2$}, Math. Res. Lett. \textbf{12} (2005), p. 897-910.

\bibitem{Ko2} M.Kol\'a\v r : {Local symmetries of finite type hypersurfaces in $\mathbb C^2$}, 
Sci. China A \textbf{49} (2006), p.
1633-1641.

\bibitem{KL} N.G.Kruzhilin, A.V.Loboda : {
Linearization of local automorphisms of pseudoconvex surfaces,}
Dokl. Akad. Nauk SSSR
 \textbf{271} (1983), p. 280-282.

\bibitem{Kow} R.Kowalski :  \textit{A hypersurface in $\Bbb C\sp 2$
whose stability group is not determined by 2-jets.}
  Proc. Amer. Math. Soc.  \textbf{130}  (2002), p. 3679--3686.

\bibitem{MZ} D.Montgomery, L.Zippin : \textit{
Topological transformation groups,}
 Interscience (1955).

\bibitem{Po} H.Poincar\'e : \textit{Les fonctions analytique de
deux variables et la repr\'esentation conforme } Rend. Circ. Mat.
Palermo \textbf{23} (1907), p. 185-220.

\bibitem{Sp} G.Schmalz, A.Spiro : \textit{ Explicit construction of a Chern-Moser
connection for CR manifolds of codimension two}  Ann. Mat. Pura Appl.
\textbf{4}
(2006), p. 337-379.

\bibitem{Sl} G.Schmalz, J.Slov\'ak : \textit{
The geometry of hyperbolic and elliptic CR-manifolds of codimension two }
Asian J. Math. \textbf{4} (2000), no. 3, p. 565-597.


\bibitem{S} N.Stanton : \textit{A normal form for rigid hypersurfaces in
$\Bbb C^2$ }, Amer. J. Math.  \textbf{113} (1991), p. 877-910.

\bibitem{V} A.G.Vitushkin : \textit{Real analytic
hypersurfaces in complex manifolds}, Russ. Math. Surv. \textbf{40}
(1985), p. 1-35.

\bibitem {W} S.M.Webster : \textit{On
the Moser normal form at a non-umbilic point}, Math. Ann.
 \textbf{233} (1978), p. 97-102.

\bibitem {Wo}  P.Wong : \textit{A construction of normal forms for weakly
pseudoconvex CR manifolds in $\Bbb C^2$ }, Invent. Math.
\textbf{69}(1982), p. 311-329.

\end{thebibliography}
\end{document}